\documentclass[twocolumn]{autart}    

\usepackage{algorithm}
\usepackage{amsmath,amssymb}
\usepackage{graphicx} 
\usepackage{xcolor}
\usepackage{cite}
\usepackage{centernot}
\usepackage{booktabs}

\let\classAND\AND
\let\AND\relax

\usepackage{algorithmic}

\let\AND\classAND
\AtBeginEnvironment{algorithmic}{\let\AND\algoAND}
	
\floatname{algorithm}{Algorithm}

\newtheorem{ass}{Assumption}
\newtheorem{lemma}{Lemma}
\newtheorem{remark}{Remark}

\DeclareMathOperator*{\argmin}{\arg\min}

\renewcommand{\epsilon}{\varepsilon}
\renewcommand{\theta}{\vartheta}
\renewcommand{\phi}{\varphi}

\newcommand{\hqed}{\hfill\qed}

\newcommand{\im}{i=1,\dots,m}
\newcommand{\sumim}{\sum_{i=1}^m}
\newcommand{\dsum}{\displaystyle \sum}
\newcommand{\T}{^\top}
\newcommand{\conv}[1]{\mathrm{conv}(#1)}

\newcommand{\vertt}[1]{\mathrm{vert}(#1)}
\newcommand{\LP}{\text{\tiny{LP}}}
\newcommand{\tunderbar}[1]{\mkern1.5mu\underline{\mkern-1.5mu#1\mkern-1.5mu}\mkern 1.5mu}

\newcommand{\Z}{\mathbb{Z}}
\newcommand{\R}{\mathbb{R}}
\newcommand{\cD}{\mathcal{D}}
\newcommand{\cP}{\mathcal{P}}
\newcommand{\cM}{\mathcal{M}}
\newcommand{\cV}{\mathcal{V}}

\newcommand{\lambdainit}{\lambda_{\mathrm{ref}}}
\newcommand{\dual}{\varphi}

\newcommand{\xfeas}{\hat{x}}
\newcommand{\slater}{\tilde{x}}
\newcommand{\ls}{\lambda^\star}
\newcommand{\Ls}{\Lambda^\star}
\newcommand{\lup}{{\overline{\lambda}}}
\newcommand{\ldn}{{\underline{\lambda}}}

\newcommand{\lsdn}{{\underline{\lambda}^\star}}
\newcommand{\iter}{k}
\newcommand{\kin}{\kappa}
\newcommand{\ko}{h}

\definecolor{mygray}{RGB}{160 160 160}

\newcommand{\forComp}{\renewcommand\algorithmicdo{\textbf{compute}}}

\newcommand{\rev}{\color{black}}
\newcommand{\btb}{\color{black}}

\begin{document}
	
	\begin{frontmatter}
		
		\title{DualBi: A dual bisection algorithm for non-convex problems with a scalar complicating constraint} 
		
		\thanks[footnoteinfo]{This paper was not presented at any IFAC 
			meeting. Corresponding author L.~Manieri. Tel. +39 02 2399 4028. }
		
		\author[polimi]{Lucrezia Manieri}\ead{lucrezia.manieri@polimi.it},    
		\author[polimi]{Alessandro Falsone}\ead{alessandro.falsone@polimi.it},               
		\author[polimi]{Maria Prandini}\ead{maria.prandini@polimi.it}  
		
		\address[polimi]{Dipartimento di Elettronica Informazione e Bioingegneria, Politecnico di Milano, Via Ponzio 34/5, 20133 Milano, Italy}  
			

		\begin{abstract}                          
			This paper addresses non-convex constrained optimization problems that are characterized by a scalar complicating constraint. We propose an iterative bisection method for the dual problem (DualBi Algorithm) that recovers a feasible primal solution, with a performance that is progressively improving throughout iterations. Application to multi-agent problems with a scalar coupling constraint results in a decentralized resolution scheme where a central unit is in charge of the update of the (scalar) dual variable while agents compute their local primal variables. In the case of multi-agent MILPs, simulations showcase the performance of the proposed method compared with state-of-the-art duality-based approaches.
		\end{abstract}
		
	\end{frontmatter}
	
\section{Introduction} 
\label{sec:introduction}
We consider non-convex constrained optimization programs that are characterized by a scalar constraint that \rev complicates their resolution.  The problem has 
the following form
\begin{align*}
	\min_x  \quad &f(x) \tag{$\cP$} \label{eq:problem} \\
	\text{subject to:} \quad &v(x) \le 0 \\
	&x \in X,
\end{align*}
where $x \in \R^n$ is the decision vector taking values in a compact (possibly non-convex) set $X \subset \R^n$, while $f: X \to \R$ and  $v: X \to \R$ are continuous scalar functions defining, respectively, the cost and the \textit{complicating} constraint. 

Continuity of $f(\cdot)$ and $v(\cdot)$ guarantees that the constrained minimum exists given the compactness of the constraint set $X \cap \{x\in \R^n: v(x) \le 0 \}$. We shall denote it as $f^\star$. We also assume that 
\begin{align}\label{eq:scalar}
    v(x) \le 0
\end{align}
is not redundant, meaning that removing it from the formulation would make the solution of the resulting problem super-optimal, i.e., $\min_{x \in X}  f(x) < f^\star$. 
\btb 

Relevant examples \rev of optimization programs in the form of~\ref{eq:problem} include resource allocation problems~\cite{katoh1998resource} with a single budget constraint, portfolio optimization~\cite{baumann_portfolio-optimization_2013}, one-dimensional cutting stocks problems~\cite{gilmore_linear_1961}, asset-backed securitization problems~\cite{mansini_securitization_2004}, and some problems in the energy domain
like e.g., economic dispatch problems~\cite{Gaing2003,manieri_dual_2024}, where multiple generators must be coordinated to meet a total power demand while accounting for their individual constraints. 
\btb

\rev The literature on non-convex optimization resolution methods that apply to~\ref{eq:problem} is quite vast. Here, we confine our discussion to those approaches that are able to preserve the feasibility of the explored candidate solutions, a feature that is crucial for many control-related applications that require ready-to-apply solutions in case of a premature interruption of the procedure. Such resolution schemes include Feasible Interior Point methods~\cite{byrd_feasible_2003}, Feasible Sequential Quadratic Programming~\cite{lawrence_computationally_2001}, and schemes based on Successive Convex Approximation~\cite{scutari_parallel_2016}.

Interior Point methods~\cite{conn_primal-dual_2000,vanderbei_interior-point_1999,wachter_implementation_2006} search for a solution of a non-convex constrained optimization program by exploring the interior of its feasible region. They typically lift the inequality constraints to the cost through a suitable barrier term and solve the resulting auxiliary barrier problem iteratively. Among these approaches, slack-based Interior Point Methods~\cite{byrd_feasible_2003} can be suitably modified to grant feasibility of all the candidate primal solutions explored throughout the iterations. Their application is, however, limited to optimization programs in which the cost and constraint functions are differentiable and may get stuck when applied to problems having both equality and inequality constraints.

Feasible Sequential Quadratic Programming (FSQP) methods~\cite{polak_computational_1971,panier_combining_1993,beck_sequential_2010}, instead, search for a solution of the problem by solving a sequence of local quadratic approximations based on first and second order information on the cost and constraints. Differently from standard Sequential Quadratic Programming approaches, these feasible variants ensure that the explored candidate primal solutions remain within the feasible set throughout the iterations by \textit{tilting}~\cite{polak_computational_1971} and \textit{bending}~\cite{panier_combining_1993,beck_sequential_2010} the search direction at each iteration. Despite the efforts to reduce computational complexity in, e.g.,~\cite{beck_sequential_2010}, FSQP methods still prove computationally intensive, since they require solving at least two QPs at each iteration.

In a similar but more recent stream of work, methods based on Successive Convex Approximation (SCA)~\cite{scutari_parallel_2016} address non-convex constrained optimization problems by solving a sequence of (simpler) sub-programs in which the non-convex cost and constraints of the original problem are replaced by suitable strongly convex approximations. This framework encompasses several resolution schemes that differ in the choice of the approximating cost and constraints, and in the algorithm used to solve the resulting convex approximation. Some examples include gradient or Newton-type methods, block coordinate descent schemes and convex-concave approximation methods , with the possibility to parallelize computations for structured problems. However, convergence is guaranteed for optimization programs with sufficiently smooth cost and constraint functions, provided that the approximate cost and constraint satisfy some suitable convexity and growth conditions. 

All these state-of-the-art approaches are devised for general non-convex problems, thus ignoring the structure of~\ref{eq:problem} and possibly resulting in an unnecessary additional computational burden.
In this work, we therefore propose an iterative scheme that is tailored to the structure of~\ref{eq:problem}. It resorts to Lagrangian duality theory to handle the complicating constraint by lifting it to the cost function, with a weight given by a scalar Lagrange multiplier. 

The strategy of lifting to the cost 
those constraints that complicate the resolution of a problem has been widely adopted in the literature. Although some theoretical results are stronger in a convex setting, there are several recent  
and also not so recent
attempts to exploit Lagrangian duality in a non-convex framework and recover tractability at least for specific problem structures. 
Examples include \cite{gasimov_augmented_2002} which proposes an approach for non-convex minimization problems with \textit{complicating} equality constraints, \cite{shi_penalty_2020} that addresses non-smooth optimization problems with partially separable structure and equality coupling constraints, and \cite{mancino-ball_decentralized_2023} which proposes a primal-dual strategy for a class of non-convex non-smooth consensus problems. 

Similarly, to~\cite{gasimov_augmented_2002,shi_penalty_2020,mancino-ball_decentralized_2023}, we propose a duality-based approach that leverages the specific structure of $\cP$ to formulate a one-dimensional dual problem. Its scalar nature is then exploited to compute an optimal dual solution by searching for a zero of the differential mapping of the dual function via a bisection method. The resulting Dual Bisection (DualBi) Algorithm is shown to either converge to an optimal primal solution or generate a sequence of feasible primal solutions with non-deteriorating performance. Although limited to problems with a single complicating constraint, the approach rests on milder assumptions on the cost and constraint functions than those in~\cite{byrd_feasible_2003,beck_sequential_2010,scutari_parallel_2016}, making it fit for 
a wider range of instances of~\ref{eq:problem} than its competitors. In addition, any existing state-of-the-art scheme for non-convex programs may be combined with the proposed approach to address the (easier but still possibly non-convex) problem obtained by lifting the complicating constraint to the cost, that needs to be solved at each iteration of the proposed procedure. \btb

The remainder of the paper is structured as follows. After a short paragraph reviewing the adopted notation, Section~\ref{sec:problem} \rev describes the proposed DualBi Algorithm and states its feasibility and performance properties, \btb whose proof is \rev deferred to \btb Appendix~\ref{app:theory}. Section~\ref{sec:multi_agent} focuses on the application of the DualBi Algorithm to multi-agent problems and discusses its advantages with respect to \rev other duality-based \btb state-of-the-art approaches,  and, in particular, to a modified version of the algorithm in~\cite{Manieri-IFAC-2023}, whose properties and resulting superiority with respect to the other methods in the literature are formally proven in Appendix~\ref{app:multi-agent}, with reference to the considered case of scalar coupling.
The performance of the proposed procedure is assessed via numerical simulations in Section~\ref{sec:numerical_analysis}.
Finally, Section~\ref{sec: conclusions} concludes the paper.

\paragraph*{Notation}
We denote by $\R$ the set of real numbers. The Cartesian product between sets is denoted with $\times$. For a function $f :~ \R^n \to \R$ we denote by $\partial f(x) \subset \R^n$ the sub-differential (i.e., the set of all sub-gradients) of $f$ at $x$. If $f$ is differentiable at $x$, then $\partial f(x) = \{ \nabla f(x) \}$ where $\nabla f(x) \in \R^n$ is the gradient of $f$ at $x$.

\section{\rev Setting  \btb and proposed solution} 
\label{sec:problem}

\rev Before delving into the proposed approach, we shall cover some preliminary facts and assumptions that are key for the adopted methodology. 

\rev As anticipated in the introduction, we consider the case when~\ref{eq:problem} is not necessarily convex, and \eqref{eq:scalar} is a complicating constraint, meaning that problem~\ref{eq:problem} is hard to solve, while the problem obtained by lifting the constraint in the cost through a penalization term is easier.\btb

In such a setting, it is natural to resort to duality theory to handle the complicating constraint. This entails introducing a (single) Lagrange multiplier $\lambda \in \R$ and define the Lagrangian function $L: \R^n \times \R \to \R$ as
\begin{equation} \label{eq:lagrangian}
	L(x,\lambda) = f(x) + \lambda v(x).
\end{equation}
We can then construct the dual function $\dual: \R \to \R$ 
\begin{equation} \label{eq:dual_function}
	\dual(\lambda) = \min_{x \in X} L(x,\lambda),
\end{equation}
which is well-defined given the continuity of functions $f$ and $v$ and the compactness of $X$, and pose the dual problem
\begin{equation*} \label{eq:dual_problem}
	\max_{\lambda \ge 0} \dual(\lambda), \tag{$\cD$}
\end{equation*}
whose optimal value (if it exists) is known to yield a lower bound for the optimal value of~\ref{eq:problem}, see, e.g.,~\cite[Section~5.1.3]{boyd2004convex}:
\begin{align}\label{eq:lower}
	 \max_{\lambda \ge 0} \btb \dual(\lambda) \le f^\star.
\end{align}
\rev To ensure that the maximum in~\ref{eq:dual_problem} is attained, we impose the following (non restrictive) assumption

\begin{ass}
\label{ass:strict-feasibility}
There exists $\tilde{x}\in X$ such that $v(\tilde x)<0$. 
\end{ass}

Under Assumption~\ref{ass:strict-feasibility}, Lemma~\ref{lemma:existence} in Appendix~\ref{app:preliminaries-A} guarantees that the set of dual optimal solution $\Ls$ is non-empty and bounded. 

\btb 

One may thus be tempted to compute an optimal solution \rev $\ls \in \Ls$ \btb of~\ref{eq:dual_problem} and then recover a primal solution of~\ref{eq:problem} as
\begin{equation} \label{eq:primal_recovery}
	x_{\lambda} \in X(\lambda) = \argmin_{x \in X} L(x,\lambda) 
\end{equation}
with $\lambda = \ls$. Unfortunately, this strategy is not guaranteed to provide an optimal primal solution or even a feasible one for the dualized constraint $v(x) \le 0$\rev, see, e.g., \cite[Example~1]{Falsone2019}.\btb

\rev In this work, \btb we show that we can exploit the dual of the original problem~\ref{eq:problem} 
to find a feasible primal solution. Starting from the observation that solving~\ref{eq:dual_problem} amounts to finding the zeros of the sub-differential of $\dual(\lambda)$, we leverage the scalar nature of the complicating constraint~\eqref{eq:scalar} and the monotonicity of the sub-differential to devise an algorithm that first finds a feasible primal solution to~\ref{eq:problem} and then solves the dual problem~\ref{eq:dual_problem} via bisection,  while preserving feasibility of the primal solution and progressively decreasing its sub-optimality level. We also show that if the set $\Ls$ of optimal dual solutions \rev is not a singleton, \btb
then an optimal primal solution is found in a finite number of steps.

\begin{algorithm}[t]
\caption{Dual Bisection (DualBi)  Algorithm
} \label{alg:dual-bisection}
\textbf{Input:}  $\lambdainit>0$, 
\label{step:lb_init}
\begin{algorithmic}[1]
	\STATE $\iter \gets 0$
	\STATE $\ldn(0) \gets 0$; 
 $\lup(0) \gets \lambdainit$; 
 \STATE $\displaystyle x_{\lup(0)} \in \argmin_{x \in X} f(x) + \lup(0) v(x)$ \hspace{-1em} \label{step:xfeas_init}
	
	\textcolor{mygray}{\%Find a feasible solution}
	\WHILE{$v(x_{\lup(k)}) > 0$} \label{step:cycle_1_start}
	\STATE $\ldn(k+1) \gets \lup(k)$  \label{step:ldn_increase} 
	\STATE $\lup(k+1) \gets 2 \lup(k)$ \label{step:lup_increase} 
	\STATE $\displaystyle x_{\lup(k+1)} \in \argmin_{x \in X} f(x) + \lup(k+1) v(x)$ \label{eq:cycle_1_find_xi}
	\STATE $\iter \gets \iter + 1$
	\ENDWHILE \label{step:cycle_1_end}
	\STATE $K \gets k$
	\\\vspace{1mm}
	\textcolor{mygray}{\%Check if solution is optimal and return it}
        \IF{$v(x_{\lup(k)}) = 0$} \label{step:check_opt_start}
	\RETURN $x_{\lup(k)}$
	\ENDIF\label{step:check_opt_end}
        \\\vspace{1mm}
		\textcolor{mygray}{\%Store feasible solution and improve it via bisection}
	\STATE $\xfeas(k) \gets x_{\lup(k)}$
	\REPEAT \label{step:cycle_2_start}
	\STATE $\lambda(\iter) \gets \frac{1}{2} \big( \lup(\iter) + \ldn(\iter) \big) $ \label{step:bisect}
	\STATE $\displaystyle x_{\lambda(\iter)} \in \argmin_{x \in X} f(x) + \lambda(\iter) v(x)$ \label{step:primal_update}
	\IF{$v(x_{\lambda(\iter)}) = 0$} \label{step:check_tight}
	\RETURN $x_{\lambda(\iter)}$ \label{step:return_tight}
	\ELSIF{$v(x_{\lambda(\iter)}) < 0$} \label{step:check_feasibility}
	\STATE $\xfeas(\iter+1) \gets x_{\lambda(\iter)}$ \label{step:xfeas_update}
	\STATE $\lup(\iter+1) \gets \lambda(\iter)$ \label{step:lup_update}
	\STATE $\ldn(\iter+1) \gets \ldn(\iter)$ \label{step:ldn_keep}
	\ELSIF{$v(x_{\lambda(\iter)}) > 0$}
	\STATE $\xfeas(\iter+1) \gets \xfeas(\iter)$ \label{step:xfeas_keep}
	\STATE $\lup(\iter+1) \gets \lup(\iter)$ \label{step:lup_keep}
	\STATE $\ldn(\iter+1) \gets \lambda(\iter)$ \label{step:ldn_update}
	\ENDIF \label{step:end_check_feasibility}
	\STATE $\iter \gets \iter + 1$
	\UNTIL{some stopping condition is satisfied.} \label{step:cycle_2_end}
\end{algorithmic}
\textbf{Output:} $\xfeas(\iter)$ 
\end{algorithm}

The proposed Algorithm~\ref{alg:dual-bisection} is structured into two cycles. The first one (cf. Steps \ref{step:cycle_1_start}-\ref{step:cycle_1_end}) seeks a feasible primal solution together with two values  $\ldn$ and $\lup$ such that $\ldn < \lup$ and the associated primal solutions satisfy $v(x_{\ldn}) > 0$ and $v(x_{\lup}) \le 0$, respectively, initializing $\ldn = 0$ and setting $\lup$ equal to a tentative value $\lambdainit >0$. If $v(x_{\lup}) > 0$ the scheme sets $\ldn$ equal to $\lup$ (cf. Step~\ref{step:ldn_increase}) and then doubles $\lup$ (cf. Step~\ref{step:lup_increase}), until a $\lup$ such that $v(x_{\lup}) \le 0$ is found. 
\rev If the obtained solution $x_{\lup}$ satisfies the complicating constraint with equality, the algorithm stops and returns it (cf. Steps~\ref{step:check_opt_start}-\ref{step:check_opt_end}).
Otherwise,  in the second cycle (cf. Steps \ref{step:cycle_2_start}-\ref{step:cycle_2_end}), the interval $[\ldn,\lup]$, with $\ldn$ and $\lup$ obtained at iteration $K$, \btb  is iteratively halved to create a sequence $\{\hat{x}(\iter)\}_{\iter\ge K}$ of feasible primal solutions of \ref{eq:problem} with non-increasing cost. 

At each iteration, the procedure computes the midpoint $\lambda$ of the current interval $[\ldn,\lup]$ (cf. Step~\ref{step:bisect}) and the associated tentative primal solution $x_{\lambda}$ (cf. Step~\ref{step:primal_update}). If $x_{\lambda}$ satisfies the complicating constraint with equality, the algorithm stops and returns $x_{\lambda}$. Conversely, the extreme points $\ldn$ and $\lup$ and the feasible solution sequence $\{\hat{x}(\iter)\}_{\iter\ge K}$ are updated based on the sign of $v(x_{\lambda})$ (cf. Steps~\ref{step:check_feasibility}-\ref{step:end_check_feasibility}). 
In particular, if $v(x_{\lambda}) < 0$ the upper-bound $\lup$ is set equal to $\lambda$, the current (feasible) primal solution $x_{\lambda}$ stored in the sequence $\{\hat{x}(\iter)\}_{\iter\ge K}$, and the lower-bound $\ldn$ remains unchanged  (cf. Steps~\ref{step:xfeas_update}-\ref{step:ldn_keep}); otherwise, if $v(x_{\lambda})>0$, then the lower-bound $\ldn$ is set equal to $\lambda$, the upper-bound $\lup$ remains unchanged, and the sequence $\{\hat{x}(\iter)\}_{\iter\ge K}$ is updated with the last feasible solution found (cf. Steps~\ref{step:xfeas_keep}-\ref{step:ldn_update}).

\rev We shall show in Proposition~\ref{prop:zero-violation} that if a solution $x_\lambda$ satisfies the complicating constraint with equality, i.e. $v(x_\lambda) = 0$, then it is also optimal. Thus, \btb the algorithm halts when a solution $x_\lambda$ with a zero violation is found \rev (cf. Steps \ref{step:check_opt_start}-\ref{step:check_opt_end} and Steps~\ref{step:check_tight}-\ref{step:return_tight}); \btb  otherwise, it keeps running until some stopping condition is met. For instance, one can stop when a maximum number of iterations for the second cycle is reached or when the decreasing length of the interval $[\ldn,\lup]$ (which provably contains $\Ls$) falls below a given threshold. 

The sequences generated by Algorithm~\ref{alg:dual-bisection} satisfy 
\begin{align*}
f(\xfeas(\iter)) \le f^\star + \lup(\iter) |v(\xfeas(\iter))|,
\end{align*}
which provides an assessment of its quality. Convergence in a finite number of iterations to a feasible (and possibly optimal) primal solution is granted by the following theorem. 

\begin{thm}[Feasibility and Performance] \label{thm:main}
Under Assumption~\ref{ass:strict-feasibility}, Algorithm~\ref{alg:dual-bisection} provides a feasible solution to the primal problem \ref{eq:problem} in a finite number of iterations $K \ge 0$.  Then, either an optimal solution is found in a finite number of iterations $\bar k \ge K$ or a sequence of feasible and cost-improving iterates $\{\xfeas(\iter)\}_{\iter\ge K}$ such that
\begin{subequations} \label{eq:guarantees}
	\begin{align}
		&v(\xfeas(\iter)) \le 0 					&&k\ge K \label{eq:feasible} \\
		&f(\xfeas(\iter+1)) \le f(\xfeas(\iter))	&&k\ge K \label{eq:decreasing_cost}
	\end{align} 
\end{subequations}
is generated.
\end{thm}

The following result shows that if a feasible solution is already available, then one can initialize Algorithm~\ref{alg:dual-bisection} in such a way as to skip the first cycle.  
\begin{thm}[Initialization] \label{thm:skip_first_cycle}
If a feasible solution $\slater \in X$ such that $v(\slater) < 0$ is known, then setting
	\begin{equation} \label{eq:slater_init}
		\lambdainit = \frac{\dual(0)-f(\slater) }{v(\slater)}
	\end{equation}
	yields $K = 0$.
\end{thm}

The proofs of both theorems are provided in Appendix~\ref{app:theory}, together with some preliminary results.

We next focus on multi-agent problems in which the complicating constraint is represented by a scalar coupling constraint, and we investigate their resolution via the DualBi Algorithm.

\rev
\section{Constraint-coupled multi-agent problems} \label{sec:multi_agent}

In this section, we focus on constraint-coupled multi-agent problems where the complicating constraint~\eqref{eq:scalar} is a scalar \textit{coupling} constraint binding decision variables of different agents. The application of the DualBi Algorithm results in a decentralized resolution scheme where a central unit is in charge of updating the (scalar) dual variable, while agents compute their local primal variables.

The decomposition of the centralized problem into multiple sub-programs of smaller size is particularly interesting for defeating the combinatorial complexity of 
\btb 
constraint-coupled  multi-agent Mixed Integer Programs (MIPs)  of the form
\begin{align*}
\min_{\xi_1,\dots,\xi_m}  \quad &\sum_{i=1}^m f_i(\xi_i) \tag{$\cM$} \label{eq:multi-agent_problem} \\
\text{subject to:} \quad &\sum_{i=1}^m v_i(\xi_i) \le 0 \\
& \xi_i \in \Xi_i\subset \R^{n_{c,i}}\times \Z^{n_{d,i}}, \quad \forall i=1,\dots,m              
\end{align*}
where $\Xi_i$ is a mixed integer compact constraint set. These \rev programs formalize those problems in which $m$ agents cooperate \btb to minimize the sum of their local cost functions $f_i$ by \rev properly \btb setting the values of $n_{c,i}$ continuous and $n_{d,i}$ discrete optimization variables collected in the decision vector $\xi_i$, for $i=1,\dots,m$, \rev while also accounting for \btb local constraints modeling the individual operational limitations, coded in the sets $\Xi_i$, $i = 1,\dots,m$, and a global constraint, coded in $\sum_{i=1}^m v_i(\xi_i) \le 0$, related to the usage of a shared resource, which couples the agents' decisions and prevents the resolution of~\ref{eq:multi-agent_problem} as $m$ separate problems.
Such a coupling constraint can be regarded as a complicating constraint, and one then just needs to define $x = [\xi_1^\top \, \cdots \, \xi_m^\top]^\top$, $X = \Xi_1 \times \cdots \times \Xi_m$, $f(x) = \sum_{i=1}^m f_i(\xi_i)$ and $v(x) = \sum_{i=1}^m v_i(\xi_i)$ to make~\ref{eq:multi-agent_problem} fit the structure of~\ref{eq:problem} and apply Algorithm \ref{alg:dual-bisection} proposed in this paper, with compactness of $X$  inherited by the compactness of the local constraint sets $\Xi_i$, $i=1,\dots,m$. 

\begin{remark}[decentralized implementation]\label{rem:decentralized}
Since the coupling constraint is the sum of the agents' local contributions and lifting it to the cost renders~\ref{eq:multi-agent_problem} separable into $m$ sub-problems, then Algorithm \ref{alg:dual-bisection} can be implemented according to a decentralized scheme where the minimization in Steps~\ref{eq:cycle_1_find_xi} and~\ref{step:primal_update} is distributed across the agents, whilst a central unit is only in charge of updating the $\lambda$-interval at each iteration.      
\end{remark}
 
As for the continuity assumption of functions $f$ and $v$ over $X$, \rev it is sufficient to require \btb continuity with respect to the $n_c=\sum_{i=1}^m n_{c,i}$ continuous variables of the decision vector $x$ over each slice of $X$ obtained by fixing the $n_d=\sum_{i=1}^m n_{d,i}$ discrete variables.

Recently,  \rev duality-based methods have been proposed in~\cite{Falsone2018,Falsone2019,LaBella2021,Manieri2023,Manieri-IFAC-2023} \btb to address constraint-coupled multi-agent MILPs (Mixed Integer Linear Programs), which are a special instance of problems in the form of~\ref{eq:multi-agent_problem} characterized by  linear cost and coupling constraint functions and with local mixed-integer polyhedral sets:
\begin{subequations}
	\label{eq:multi-agent_MILP}
	\begin{align}
		\label{eq:multi-agent MILP cost function} \min_{\xi_1,\dots,\xi_m} \quad &\sum_{i=1}^{m} c_i^\top \xi_i \\
		\label{eq:multi-agent MILP coupling constraints} \text{subject to:} \quad &\sum_{i=1}^{m} A_i \xi_i -b \leq 0	\\
		\label{eq:multi-agent MILP local constraints} &\xi_i \in \Xi_i, \quad i=1,\dots,m,
	\end{align}
\end{subequations}
where $\Xi_i=\{\xi_i \in \R^{n_{c,i}}\times\Z^{n_{d,i}} : D_i x_i \le d_i\}$, with matrix $D_i$ and vector $d_i$ of appropriate dimensions. \rev Despite their linearity, these problems can capture a wide range of engineering applications, including optimal operation of hybrid systems with linear continuous dynamics, finite-state machines, bilinear or PieceWise Affine (PWA) dynamical systems subject to PWA/linear constraints and any other system that can be described as a multi-agent Mixed Logical Dynamical~\cite{bemporad_control_1999} system, where each unit is characterized by interleaved physical and digital components. As a result, multi-agent MILPs in the form of \eqref{eq:multi-agent_MILP} arise in several domains ranging from energy systems \cite{LaBella2021} and building management \cite{ioli2015iterative} to multi-robot vehicle coordination \cite{camisa2022multi}. \btb

\rev State-of-the-art resolution schemes for \eqref{eq:multi-agent_MILP} resort to a tightening of the coupling constraint to recover a feasible primal solution, which typically results in a worse performance than the DualBi Algorithm proposed in this work, which does not resort to any tightening of the coupling constraint to enforce feasibility.  
More specifically, in those works, \btb an iterative decentralized scheme is proposed  that computes a feasible primal solution  by solving the dual of the following modified problem
\begin{subequations}
	\label{eq:multi-agent_tightprimal}
	\begin{align}
		\label{eq:multi-agent_tightprimal-cost} \min_{\xi_1,\dots,\xi_m} \quad &\sum_{i=1}^{m} c_i^\top \xi_i \\
		\label{eq:multi-agent_tightprimal-coupling-constraints} \text{subject to:} \quad &\sum_{i=1}^{m} A_i \xi_i -b \leq -\rho	\\
		\label{eq:multi-agent_tightprimal-local-constraints} &\xi_i \in \Xi_i, \quad i=1,\dots,m,
	\end{align}
\end{subequations}
where the coupling constraint is fictitiously tightened of a quantity $\rho \ge 0$. In \cite{Falsone2019}, the tightening vector $\rho$ is adaptively increased at each iteration in a decentralized fashion based on the previously explored solutions. Convergence to a feasible primal solution in a finite number of iterations is guaranteed if the feasibility sets $\Xi_i$ are non-empty and bounded for $i=1,\dots,m$ and if a given relaxation of the primal and dual tightened problems admit a unique solution when $\rho$ is set equal to its convergence value $\bar{\rho}$. A bound on the sub-optimality level is also provided and is shown to worsen as the infinity norm of $\bar{\rho}$ increases.
In \cite{Falsone2018}, the same approach is implemented in a distributed fashion while preserving the feasibility guarantees and the bound on the quality of the retrieved primal solution. In \cite{LaBella2021}, instead, the approach in \cite{Falsone2019} is made less conservative by making the increase of $\rho$ more cautious, but no guarantees are provided. A variant of the method in \cite{LaBella2021} is proposed in \cite{Manieri2023} where similar guarantees as in  \cite{Falsone2019} are given but with improved performance results. 

In \cite{Manieri-IFAC-2023}, conservativeness is further reduced by allowing a decrease of the tightening across the iterations. The proposed procedure (cf. \cite[Algorithm~1]{Manieri-IFAC-2023}) updates the tightening vector $\rho$ based on the mismatch between a mixed-integer solution satisfying
\begin{equation}
\label{eq:xi_lsr}
    \xi_{\ls_\rho} \in \argmin_{\Xi_1 \times \cdots \Xi_m} \left( [c_1\T \cdots c_m\T]+{\ls_\rho}\T [A_1 \cdots A_m]\right) \xi \, ,
\end{equation}
$\ls_\rho$ being the ($\rho$-dependent) optimal solution of the dual of~\eqref{eq:multi-agent_tightprimal}, and the solution $\xi^\LP_\rho$ of the following relaxation of~\eqref{eq:multi-agent_tightprimal}  
\begin{subequations}
	\label{eq:multi-agent_relaxed_tightprimal}
	\begin{align}
		\label{eq:multi-agent-rel-tp-cost} \min_{\xi_1,\dots,\xi_m} \quad &\sum_{i=1}^{m} c_i^\top \xi_i \\
		\label{eq:multi-agent-rel-tp-coupling-constraints} \text{subject to:} \quad &\sum_{i=1}^{m} A_i \xi_i -b \leq -\rho	\\
		\label{eq:multi-agent-rel-tp-local-constraints} &\xi_i \in \conv{\Xi_i}, \quad i=1,\dots,m,
	\end{align}
\end{subequations}
which is obtained by replacing the local sets $\Xi_i,\, i=1,\dots,m$ with their convex hulls and shares the same dual problem with~\eqref{eq:multi-agent_tightprimal}.

In~\cite{Manieri-IFAC-2023}, a separate algorithm (cf.~\cite[Algorithm~2]{Manieri-IFAC-2023}) is proposed to compute the pair of solutions $(\xi_{\ls_\rho},\xi^\LP_\rho)$ used to update $\rho$. We shall show in Proposition~\ref{prop:ifac-inner-loop-scalar} (whose proof is in Appendix~\ref{app:multi-agent}) that if the considered multi-agent MILP has a single coupling constraint, then ~\cite[Algorithm~2]{Manieri-IFAC-2023} explores feasible solutions even when no tightening is introduced ($\rho=0$). We \rev thus \btb propose in Algorithm~\ref{algo:ifac-inner} a modified version of~\cite[Algorithm~2]{Manieri-IFAC-2023} for the case of a scalar coupling constraint of interest, which differs only for the rule used to select $\xi_{\ls_\rho}$, so that one of these feasible solutions is returned to~\cite[Algorithm~1]{Manieri-IFAC-2023} and the procedure can be halted after the first iteration since $\rho$ will keep being equal to zero indefinitely.

\rev In fact, Algorithm~\ref{algo:ifac-inner} is \btb identical to~\cite[Algorithm~2]{Manieri-IFAC-2023} except for Steps~\ref{step:ifac-tau-selection} and~\ref{step:ifac-inner-sel-feas-mixed}
where the mixed-integer solution $\xi_{\ls_\rho}$ to be returned is selected. \rev We briefly describe it next to make the paper self-contained. Steps~\ref{step:begin_subg}-\ref{step:end_subg} apply the well-established sub-gradient method to solve the dual of~\eqref{eq:multi-agent_tightprimal} (and~\eqref{eq:multi-agent_relaxed_tightprimal}) for the current value of $\rho$ by alternating a primal and a dual update (cf. Step~\ref{step:XLP_tentativePrimal} and \ref{step:end_subg}, respectively) up to convergence of the dual variable $\lambda$, which is guaranteed 
when the step-size $\alpha(\kappa)$ in Step~\ref{step:end_subg} satisfies 
 \begin{align}\label{eq:alphak}
 \sum_{\kappa=0}^\infty \alpha(\kappa)=\infty && \text{and} && \sum_{\kappa=0}^\infty \alpha^2(\kappa)<\infty
 \end{align}
(see, e.g., \cite[Theorem 2.2]{shor2012minimization} for a proof). The primal solutions computed in Step~\ref{step:XLP_tentativePrimal}, on the other hand, are not guaranteed to converge to an optimal solution $\xi_\rho^\LP$ of~\eqref{eq:multi-agent_relaxed_tightprimal}. For this reason, Steps~\ref{step:ifac-inner-xi-hat-start}-\ref{step:ifac-endrecovery} implement a recovery procedure to compute  $\xi_\rho^\LP$ by leveraging the results in~\cite[p.~117-118]{shor2012minimization} that guarantee convergence to $\xi_\rho^\LP$ by letting $h$ and $\kappa$ in Step~\ref{step:ifac-inner-xi-hat} go to $+\infty$. Finally, Steps~\ref{step:ifac-tau-selection} and~\ref{step:ifac-inner-sel-feas-mixed} select the feasible mixed-integer solution $\xi_{\ls_\rho}$ attaining the minimum cost, which is returned alongside $\xi_\rho^\LP$. 
The following result guarantees that the index $\tau^\star$ in Step~\ref{step:ifac-tau-selection}, and thus the mixed-integer solution $\xi_{\ls_\rho}$, is well-defined.

\begin{prop}
\label{prop:ifac-inner-loop-scalar}
    Consider Algorithm~\ref{algo:ifac-inner} and fix  $\rho=0$. If the sets $\Xi_i, i=1,\dots,m$ in~\eqref{eq:multi-agent MILP local constraints} are non-empty and bounded and~\eqref{eq:multi-agent MILP coupling constraints} is a scalar non-redundant constraint, then for any $\breve{\kappa}>0$ there exists $\bar{\kappa}>\breve{\kappa}$ such that the solutions $\xi_1(\bar\kappa),\dots,\xi_m(\bar\kappa)$ computed in Step~\ref{step:XLP_tentativePrimal} of Algorithm~\ref{algo:ifac-inner} satisfy
    \begin{align}
        & \sum_{i=1}^m A_i \xi_i (\bar\kappa) \le b . \label{eq:proof-scalar-eq2}
    \end{align}
 \end{prop}

\btb

\begin{remark}
We warn the attentive reader that a different notation with respect to this paper was adopted in~\cite[Algorithm~2]{Manieri-IFAC-2023}. In particular, the vector of primal decision variables and its sub-components were denoted as $x$ and $x_i$, respectively, instead of $\xi$ and $\xi_i$. 
A primal solution satisfying~\eqref{eq:xi_lsr}, was denoted with the symbol $x(\ls_\rho)$ in place of $\xi_{\ls_\rho}$. 
\end{remark}

\begin{algorithm}[t]
\caption{Strategy to compute $\xi^{\LP}_\rho$ and $\xi_{\lambda^\star_{\rho}}$ (cf.~\cite[Algorithm~2]{Manieri-IFAC-2023}}
\label{algo:ifac-inner}
\textbf{Input:} $\rho$
\begin{algorithmic}[1]
    \STATE $\ko = 1$
    \FOR{$\kin=1,2,\dots$}  \label{step:XLP_repeat}\forComp
        \FOR{$i=1,\dots,m$}  \label{step:begin_subg}
        \vspace{1mm}
            \STATE $\displaystyle \xi_i (\kin) = \xi_{i,\lambda(\kin)} \in \argmin_{\xi_i\in \Xi_i} (c_i^\top  +\lambda(\kin)^\top A_i )  \xi_i$ \label{step:XLP_tentativePrimal}
            \vspace{-1mm}
        \ENDFOR \label{step:ifac-inner-end-subg-primal}
        \vspace{1mm}
        \STATE $\mu(\kin) = \sum_{i=1}^m A_i \xi_i(\kin)-b$ \label{step:XLP_computeSubg}
        \vspace{1mm}
        \STATE $\lambda(\kin+1) = [ \lambda(\kin)+\alpha(\kin)( \mu(\kin)+\rho)]_+$ \label{step:end_subg} \label{step:ifac-dual-update}
        \vspace{1.5mm}
    \IF{ $\lambda$ at convergence for more than $w$ iterations }
        \FOR{$i=1,\dots,m$}  \label{step:ifac-inner-xi-hat-start}
            \STATE $\tilde{\xi}_i(\ko)={\dsum_{\tau=\kin-\ko}^{\kin} \alpha(\tau)\,\xi_i(\tau)}/{\dsum_{\tau=\kin-\ko}^{\kin} \alpha(\tau)}$ \label{step:ifac-inner-xi-hat}
        \ENDFOR \label{step:ifac-inner-xi-hat-end}
        \vspace{1mm}
        \STATE $\tilde{\xi}(h)=[\tilde{\xi}_1\T(h)\,\cdots\,\tilde{\xi}_m\T(h)]\T$ \label{step:ifac-endrecovery}
        \vspace{1mm}
    \IF {$\tilde{\xi}$ at convergence for $w$ iterations} 
    \vspace{1.5mm}
    \STATE ${\xi}^{\LP}_\rho = \tilde{\xi}(\kin)$ \label{step:XLP_defineXLP}
    \STATE{$\displaystyle\tau^\star \in\argmin_{\{\tau \in \mathbb{N}:\, \tau\ge \kappa-h\}} \sum_{i=1}^m c_i\T \xi_i(\tau)$} \label{step:ifac-tau-selection}\\
    \hspace{2cm} $\displaystyle \text{s.t.}\,\sumim A_i \xi_i - b \le 0$
    \vspace{2mm}
    \STATE $\xi_{\ls_\rho} = [\xi_1\T(\tau^\star) \cdots \,\xi_m\T(\tau^\star)]\T$ \label{step:ifac-inner-sel-feas-mixed}
    \vspace{1.5mm}
    \RETURN ${\xi}^{\text{\tiny LP}}_\rho, \,\xi_{\lambda^\star_\rho}$ \label{step:xLP_return}
    \ENDIF
    \STATE $\ko = \ko+1$
    \ENDIF
   \ENDFOR
\end{algorithmic}
\textbf{Output:} ${\xi}^{\text{\tiny LP}}_\rho, \,\xi_{\lambda^\star_\rho}$
\end{algorithm}

With the proposed modification, the decentralized scheme in~\cite{Manieri-IFAC-2023} computes a feasible solution of~\eqref{eq:multi-agent_MILP} without introducing any tightening and can be halted after the first iteration. To see this, consider the complementary slackness conditions (c.f.~\cite[Section 3.4, p.301]{Bertsekas1997nonlinear}) for problem~\eqref{eq:multi-agent_relaxed_tightprimal} with $\rho=\rho(0)=0$, given by
\begin{equation}
\label{eq:slackness}
    \ls_{0} \left(\sumim A_i \,{\xi^{\LP}_{i,{0}}}  - b\right) =0 .
\end{equation}
Under the standing assumption that~\eqref{eq:multi-agent MILP coupling constraints} is a \rev non-redundant \btb constraint, $\ls_0$ satisfies $\ls_0>0$ and, thus, by~\eqref{eq:slackness} we have that $\sumim A_i \,{\xi^{\LP}_{i,{0}}}  = b$, meaning that constraint~\eqref{eq:multi-agent-rel-tp-coupling-constraints} is active at the optimum. If we now plug this expression in the $\rho$-update in~\cite[Step~13 of Algorithm~1]{Manieri-IFAC-2023} we get 
\begin{equation*}
\label{eq:proof-ifac-rho-update-feas}
    \rho(1) =\hspace{-0.5mm}\Bigg[\sum_{i=1}^m A_i \Big(\xi_{i,\ls_0} -
    \xi^{\LP}_{i,{0}}\Big) \Bigg]_+\hspace{-2.3mm}  = \left[\dsum_{i=1}^m A_i \xi_{i,\ls_0} -  b \right]_+ \, .
\end{equation*}  
As a result, owing to Proposition~\ref{prop:ifac-inner-loop-scalar}
if the solutions $\xi_{1,\ls_0}, \,\dots \, ,\xi_{m,\ls_0}$ are computed via Algorithm~\ref{algo:ifac-inner} instead of~\cite[Algorithm~2]{Manieri-IFAC-2023}, we have that $\sum_{i=1}^m A_i\xi_{i,\ls_0}-b\le0$ and, thus, $\rho(1)=0=\rho(0)$. Therefore, the procedure can be halted since the second iteration would be equal to the first one.

In the next section,  we apply the DualBi Algorithm on some randomly generated instances of constraint-coupled multi-agent MILPs, and since the scheme in~\cite{Manieri-IFAC-2023} was already shown to outperform its competitors, we compare Algorithm~\ref{alg:dual-bisection} only with the improved version of~\cite{Manieri-IFAC-2023} proposed in this work.

\section{Numerical examples} 
\label{sec:numerical_analysis}

We consider problems in the form of the constraint-coupled multi-agent MILP  \eqref{eq:multi-agent_MILP}, where each agent has $n_{c,i}=5$ continuous variables and $n_{d,i}=3$ discrete variables. 

Vector $c_i\in\R^8$ defining the local cost function in \eqref{eq:multi-agent MILP cost function} has elements selected randomly in the interval $[-1,0]$. 
Each local mixed integer set $\Xi_i$ in \eqref{eq:multi-agent MILP local constraints} is given by 
\[
\Xi_i=\{\xi_i \in \R^{n_{c,i}}\times\Z^{n_{d,i}} : |\xi_i|\le 10\cdot\mathbf{1} \, \wedge \, G_i \xi_i \le g_i\}
\]
where $\mathbf{1}\in\R^8$ is a vector of ones, $G_i\in\R^{10\times 8}$ has entries drawn from the standard normal distribution, and $g_i\in\R^8$ has elements selected randomly in the interval $[0,1]$. 
As for the scalar coupling constraint, vectors $A_i$ in \eqref{eq:multi-agent MILP coupling constraints} have random entries in $[0,1]$, and the budget $b$ is set equal to half of the agents' total resource usage in the absence of coupling, so as to ensure that the resulting coupling constraint is \rev not redundant.\btb 

We compare the decentralized implementation of DualBi (see Remark \ref{rem:decentralized}) with the modified decentralized scheme in~\cite{Manieri-IFAC-2023}, which is the least conservative duality-based approach among \cite{Falsone2018,Falsone2019,LaBella2021,Manieri2023,Manieri-IFAC-2023} and is guaranteed to compute a feasible solution, due to the scalar nature of the coupling. 
Both algorithms are coded in MATLAB R2020b, with local MILPs solved using CPLEX v12.10. Simulations are performed on a laptop equipped with an Intel Core {i7-9750HF} CPU @2.60GHz and 16GB of RAM. 

Since $x=\mathbf{0}$ is a feasible solution for every instance of the problem, we  
set $\lambda_{\text{ref}}$ in Algorithm~\ref{alg:dual-bisection} according to \eqref{eq:slater_init} with $\tilde{x}=\mathbf{0}$, thus getting $K=0$. 
To ensure a fair comparison, we also initialize the procedure in~\cite{Manieri-IFAC-2023} by selecting $\lambda(0)=\lambdainit$. 

\begin{figure}
    \includegraphics[width=0.9\linewidth]{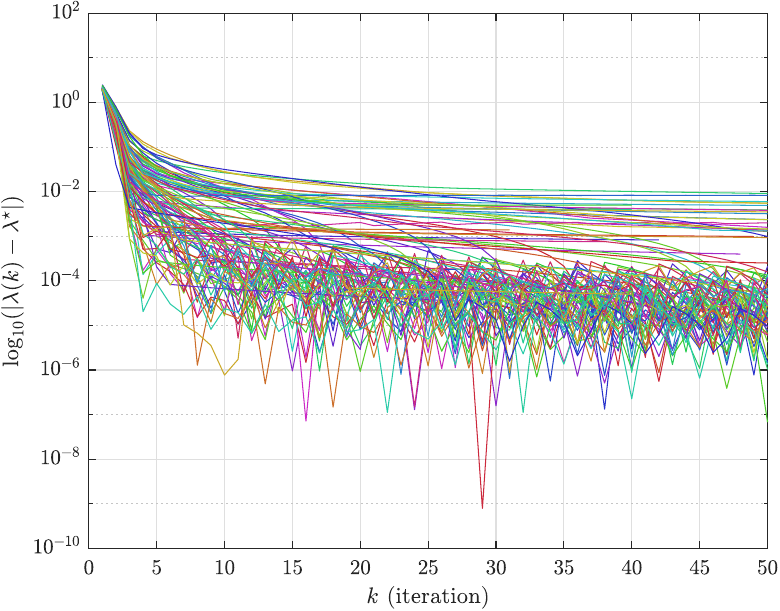}
    \caption{Absolute distance between the dual variable $\lambda(k)$ computed by the scheme in~\cite{Manieri-IFAC-2023} and $\ls$ for $k$ ranging from 0 to 50 over 100 problem instances with $m=100$ agents.}
\label{fig:test-lambda-ev}
 \end{figure}
 
We start by considering 100 problems with $m=100$ agents and 100 different parameter sets. For each test, we run Algorithm~\ref{alg:dual-bisection} until the length of the interval $[\ldn,\lup]$ falls below $10^{-5}$,  whilst the procedure in~\cite{Manieri-IFAC-2023} is executed up to convergence of the internal cycle. Specifically, we assume that practical convergence is achieved whenever the absolute difference between two consecutive values of the dual variable keeps below  $10^{-5}$ for $w=10$ consecutive iterations. As for the step-size in Step~\ref{step:ifac-dual-update} in Algorithm~\ref{algo:ifac-inner}, we set $\alpha(k)=\bar{\alpha}/k$ and we select a value $\bar{\alpha}=7\cdot10^{-4}$ such that in all tests $\lambda(k)$ decreases with a sufficiently high rate and does not exhibit oscillations for small $k$'s, as it is typically the case for excessively large step-sizes. Figure~\ref{fig:test-lambda-ev} shows the evolution of the quantity $|\lambda(k)-\ls|$ throughout the first 50 iterations of the scheme in~\cite{Manieri-IFAC-2023} for this choice of $\alpha(k)$, in all 100 tests. Since all curves approach zero with considerably small or no oscillations and the absolute distance $|\lambda(k)-\ls|$ falls below $10^{-2}$ in 50 iterations, we can conclude that the choice of $\alpha(k)$ is adequate for all problem instances.

We assess the quality of a feasible primal solution through the following index
\begin{equation}
    \label{eq:gap} \Delta f_\% = \frac{\hat{f}-\varphi^\star}{|\varphi^\star|}\, ,
\end{equation}
which measures the relative distance between the cost $\hat{f}$ attained by the feasible primal solution returned by either algorithm and the optimal dual cost $\varphi^\star$, used as a proxy for the optimal primal cost $f^\star$, due to~\eqref{eq:lower}.

Since the two procedures differ for the dual update step, and this has a negligible cost with respect to the (common) decentralized primal update, we assess the computational effort required by each scheme based on the number of iterations performed to reach the same accuracy in the solution of the dual problem, given by the difference of two consequent $\lambda$ iterates in Algorithm~\cite{Manieri-IFAC-2023} and the $\lambda$ interval in the DualBi Algorithm.

Across the 100 tests, the difference between the best cost computed by DualBi and the approach in~\cite{Manieri-IFAC-2023} never exceeds $9.10\cdot{10}^{-13}$. The average gap $\Delta f_\%$ (which is the same for the two schemes) is equal to $1.01\%$, suggesting that the primal solutions obtained by both procedures are close to the optimal ones. However, as expected, DualBi is much faster, requiring only $18$ to $19$ iterations to reach the accuracy of $10^{-5}$, whereas the scheme in~\cite{Manieri-IFAC-2023} performs approximately $300$ 
iterations on average and never less than $27$ iterations.

We now let the number of agents vary between $100$ and $1000$ to evaluate the performance of the two approaches as the size of the problem increases,  considering a single instance for each value of $m$. We use the same initialization and stopping criteria, but since the cost and the violation of the coupling constraint grow linearly with $m$, we multiply the constant $\bar{\alpha}$ in the step-size $\alpha(k)$ by a factor $m/100$.

\begin{table}
    \centering
    \begin{tabular}{cccc}
        \toprule
        $\mathbf{m}$& $\mathbf{\Delta f{\%}}$& $\mathbf{K}_D$ & $\mathbf{K}_M$  \\
        \midrule
         \textbf{100}& $1.14\cdot 10^{-1}$&18&318\\
         \textbf{250}&$2.00\cdot 10^{-3}$&18&36\\
         \textbf{500}&$1.15\cdot 10^{-4}$&18&55\\
         \textbf{750}&$3.53\cdot 10^{-2}$&18&664\\
         \textbf{1000}&$3.21\cdot 10^{-2}$&18&807\\
         \bottomrule
    \end{tabular}
      \caption{Performance of  DualBi and its competitor 
    for problems with increasing number of agents $m$ \rev in terms of relative optimality gap $\Delta f_{\%}$ and total number of iterations performed by the DualBi Algorithm and the proposed variant of the scheme in~\cite{Manieri-IFAC-2023}, denoted respectively as $K_D$ and $K_M$\btb. }
    \label{tab:perf-inc_agents}
\end{table}

Table~\ref{tab:perf-inc_agents} reports the performance obtained for $m\in\{100,250,500,750,1000\}$ in terms of $\Delta f_\%$ and the total number of iterations performed by DualBi and its competitor, respectively denoted as $K_D$ and $K_M$. Since the solutions computed by the two algorithms attain similar costs in all instances, we report a single value of $\Delta f_{\%}$ for each test. Irrespective of $m$, DualBi performs the same number of iterations and consistently less than its competitor. The computational advantage is most appreciable for the instance with $m=1000$ agents, where the scheme in~\cite{Manieri-IFAC-2023} requires a number of iterations which is 40 times that of DualBi to obtain an equally good primal solution.

As a final remark, it is worth noticing that a decentralized implementation of DualBi provides close-to-optimal solutions in less than 0.52 seconds even when the number of agents is $m=1000$, whereas for $m=100$ agents solving problem \eqref{eq:multi-agent_MILP} 
required more than 5 hours with the available computational resources.

\section{Conclusions}
\label{sec: conclusions}
This paper proposes a dual bisection (DualBi) algorithm for addressing non-convex problems with a scalar complicating constraint. 
The algorithm computes a sequence of primal iterates that are feasible and with a non-increasing cost. 
Interestingly, its application to constraint-coupled multi-agent MILPs provides a scalable decentralized resolution scheme that outperforms competing algorithms.
The main limitation of the DualBi Algorithm is that it can be applied only to the case of a scalar complicating constraint, \rev due to the nature of the bisection method \btb. The extension to the case of multiple constraints is indeed an open but quite challenging task. 

\begin{ack}                             
This paper is supported partly by the PNRR-PE-AI FAIR project funded by the NextGeneration EU program, and partly by the Italian Ministry of Enterprises and Made in Italy in the framework of the project 4DDS (4D Drone Swarms) under grant no. F/310097/01-04/X56.
\end{ack}
	
\bibliographystyle{ieeetr}        
\bibliography{bibliography}           

\begin{thebibliography}{10}

\bibitem{katoh1998resource}
N.~Katoh and T.~Ibaraki, ``Resource allocation problems,'' {\em Handbook of
  Combinatorial Optimization: Volume1--3}, pp.~905--1006, 1998.

\bibitem{baumann_portfolio-optimization_2013}
P.~Baumann and N.~Trautmann, ``Portfolio-optimization models for small
  investors,'' {\em Mathematical Methods of Operations Research}, vol.~77,
  pp.~345--356, June 2013.

\bibitem{gilmore_linear_1961}
P.~C. Gilmore and R.~E. Gomory, ``A {Linear} {Programming} {Approach} to the
  {Cutting}-{Stock} {Problem},'' {\em Operations Research}, vol.~9,
  pp.~849--859, Dec. 1961.

\bibitem{mansini_securitization_2004}
R.~Mansini and U.~Pferschy, ``{Securitization of Financial Assets:
  Approximation in Theory and Practice},'' {\em Computational Optimization and
  Applications}, vol.~29, no.~2, pp.~147--171, 2004.

\bibitem{Gaing2003}
Z.-L. Gaing, ``Particle swarm optimization to solving the economic dispatch
  considering the generator constraints,'' {\em IEEE Transactions on Power
  Systems}, vol.~18, no.~3, pp.~1187--1195, 2003.

\bibitem{manieri_dual_2024}
L.~Manieri, A.~Falsone, and M.~Prandini, ``A dual bisection approach to
  economic dispatch of generators with prohibited operating zones,'' {\em IEEE
  Control Systems Letters}, 2024.
\newblock Publisher: IEEE.

\bibitem{byrd_feasible_2003}
R.~H. Byrd, J.~Nocedal, and R.~A. Waltz, ``Feasible interior methods using
  slacks for nonlinear optimization,'' {\em Computational Optimization and
  Applications}, vol.~26, pp.~35--61, 2003.

\bibitem{lawrence_computationally_2001}
C.~T. Lawrence and A.~L. Tits, ``A {Computationally} {Efficient} {Feasible}
  {Sequential} {Quadratic} {Programming} {Algorithm},'' {\em SIAM Journal on
  Optimization}, vol.~11, pp.~1092--1118, Jan. 2001.

\bibitem{scutari_parallel_2016}
G.~Scutari, F.~Facchinei, and L.~Lampariello, ``Parallel and distributed
  methods for constrained nonconvex optimization—{Part} {I}: {Theory},'' {\em
  IEEE Transactions on Signal Processing}, vol.~65, no.~8, pp.~1929--1944,
  2016.
\newblock Publisher: IEEE.

\bibitem{conn_primal-dual_2000}
A.~R. Conn, N.~I. Gould, D.~Orban, and P.~L. Toint, ``A primal-dual
  trust-region algorithm for non-convex nonlinear programming,'' {\em
  Mathematical programming}, vol.~87, pp.~215--249, 2000.
\newblock Publisher: Springer.

\bibitem{vanderbei_interior-point_1999}
R.~J. Vanderbei and D.~F. Shanno, ``An interior-point algorithm for nonconvex
  nonlinear programming,'' {\em Computational Optimization and Applications},
  vol.~13, pp.~231--252, 1999.
\newblock Publisher: Springer.

\bibitem{wachter_implementation_2006}
A.~Wächter and L.~T. Biegler, ``On the implementation of an interior-point
  filter line-search algorithm for large-scale nonlinear programming,'' {\em
  Mathematical Programming}, vol.~106, pp.~25--57, Mar. 2006.

\bibitem{polak_computational_1971}
E.~Polak, {\em Computational methods in optimization: a unified approach},
  vol.~77.
\newblock Academic press, 1971.

\bibitem{panier_combining_1993}
E.~R. Panier and A.~L. Tits, ``On combining feasibility, descent and
  superlinear convergence in inequality constrained optimization,'' {\em
  Mathematical programming}, vol.~59, no.~1, pp.~261--276, 1993.
\newblock Publisher: Springer.

\bibitem{beck_sequential_2010}
A.~Beck, A.~Ben-Tal, and L.~Tetruashvili, ``A sequential parametric convex
  approximation method with applications to nonconvex truss topology design
  problems,'' {\em Journal of Global Optimization}, vol.~47, pp.~29--51, 2010.
\newblock Publisher: Springer.

\bibitem{gasimov_augmented_2002}
R.~N. Gasimov, ``Augmented {Lagrangian} {Duality} and {Nondifferentiable}
  {Optimization} {Methods} in {Nonconvex} {Programming},'' {\em Journal of
  Global Optimization}, vol.~24, no.~2, pp.~187--203, 2002.

\bibitem{shi_penalty_2020}
Q.~Shi and M.~Hong, ``Penalty dual decomposition method for nonsmooth nonconvex
  optimization—{Part} {I}: {Algorithms} and convergence analysis,'' {\em IEEE
  Transactions on Signal Processing}, vol.~68, pp.~4108--4122, 2020.
\newblock Publisher: IEEE.

\bibitem{mancino-ball_decentralized_2023}
G.~Mancino-Ball, Y.~Xu, and J.~Chen, ``A decentralized primal-dual framework
  for non-convex smooth consensus optimization,'' {\em IEEE Transactions on
  Signal Processing}, vol.~71, pp.~525--538, 2023.
\newblock Publisher: IEEE.

\bibitem{Manieri-IFAC-2023}
L.~Manieri, A.~Falsone, and M.~Prandini, ``A novel decentralized approach to
  large-scale multi-agent {MILP}s,'' in {\em 22nd World Congress of the
  International Federation of Automatic Control (IFAC 2023)}, 2023.
\newblock To appear on IFAC PoL.

\bibitem{boyd2004convex}
S.~Boyd, S.~P. Boyd, and L.~Vandenberghe, {\em Convex optimization}.
\newblock Cambridge university press, 2004.

\bibitem{Falsone2019}
A.~Falsone, K.~Margellos, and M.~Prandini, ``A decentralized approach to
  multi-agent {MILP}s: Finite-time feasibility and performance guarantees,''
  {\em Automatica}, vol.~103, pp.~141--150, 2019.

\bibitem{Falsone2018}
A.~Falsone, K.~Margellos, and M.~Prandini, ``A distributed iterative algorithm
  for multi-agent {MILP}s: Finite-time feasibility and performance
  characterization,'' {\em IEEE Control Systems Letters}, vol.~2, no.~4,
  pp.~563--568, 2018.

\bibitem{LaBella2021}
A.~{La Bella}, A.~Falsone, D.~Ioli, M.~Prandini, and R.~Scattolini, ``A
  mixed-integer distributed approach to prosumers aggregation for providing
  balancing services,'' {\em International Journal of Electrical Power {\&}
  Energy Systems}, vol.~133, p.~107228, dec 2021.

\bibitem{Manieri2023}
L.~Manieri, A.~Falsone, and M.~Prandini, ``Handling complexity in large scale
  cyber-physical systems through distributed computation,'' 2023.
\newblock To appear as a chapter of a book.

\bibitem{bemporad_control_1999}
A.~Bemporad and M.~Morari, ``Control of systems integrating logic, dynamics,
  and constraints,'' {\em Automatica}, vol.~35, no.~3, pp.~407--427, 1999.
\newblock Publisher: Elsevier.

\bibitem{ioli2015iterative}
D.~Ioli, A.~Falsone, and M.~Prandini, ``An iterative scheme to hierarchically
  structured optimal energy management of a microgrid,'' in {\em 2015 54th IEEE
  Conference on Decision and Control (CDC)}, pp.~5227--5232, IEEE, 2015.

\bibitem{camisa2022multi}
A.~Camisa, A.~Testa, and G.~Notarstefano, ``Multi-robot pickup and delivery via
  distributed resource allocation,'' {\em IEEE Transactions on Robotics},
  vol.~39, no.~2, pp.~1106--1118, 2023.

\bibitem{shor2012minimization}
N.~Shor, K.~C. Kiwiel, and A.~Ruszcaynski, ``{Minimization methods for
  non-differentiable functions},'' 1985.

\bibitem{Bertsekas1997nonlinear}
D.~P. Bertsekas, ``Nonlinear programming,'' {\em Journal of the Operational
  Research Society}, vol.~48, no.~3, pp.~334--334, 1997.

\bibitem{rockafellar1970convex}
R.~T. Rockafellar, {\em Convex analysis}, vol.~18.
\newblock Princeton university press, 1970.

\end{thebibliography}
\appendix

\section{Theoretical analysis of the DualBi Algorithm} \label{app:theory}
\subsection{Preliminaries} \label{app:preliminaries-A}
Let us first note that the dual function $\dual(\lambda)$ is the point-wise minimum of a collection of affine functions (cf.~\eqref{eq:lagrangian} and~\eqref{eq:dual_function}) of $\lambda$ and, as such, is concave in $\lambda$ even when the primal problem is non-convex, see, e.g., \cite[Section~5.1.2]{boyd2004convex}. This entails that every local maximum is also a global maximum.

\rev
We start by proving that the maximum in~\ref{eq:dual_problem} is actually attained under Assumption \ref{ass:strict-feasibility}. 
\begin{lemma}    \label{lemma:existence}
Under Assumption~\ref{ass:strict-feasibility}, the set $\Ls$ of the optimal solutions of~\ref{eq:dual_problem} is non-empty and bounded. 
    \begin{pf}
By definition \eqref{eq:primal_recovery}, $x_{\lambda}$ satisfies 
\begin{equation*}
    f(x_\lambda) +\lambda v(x_\lambda) \le f(x)+
    \lambda v(x) \qquad \forall x \in X,
\end{equation*}
which can be re-arranged as 
\begin{equation}
\label{eq:lambda_ineq}
v(x_\lambda)-v(x)  \le \frac{f(x)
-f(x_\lambda)}{\lambda}\le \frac{\gamma}{\lambda}, 
\end{equation}
for all $x \in X$ and $\lambda>0$, with 
\[
\gamma= \max_{x\in X} f(x)- \min_{x\in X}f(x)<\infty\, ,\]
since $f(\cdot)$ is continuous and $X$ is compact.  
If we let $\lambda\to +\infty$ in \eqref{eq:lambda_ineq},  then we get 
\begin{equation*}
    \limsup_{\lambda\to+\infty} v(x_\lambda) \le v(x) \quad \forall x\in X,  
\end{equation*} 
which implies that 
$ \limsup_{\lambda\to+\infty} v(x_\lambda) \le v(\tilde x) <0$,   
$\tilde x$ being defined in Assumption \ref{ass:strict-feasibility}.  
This entails that there exists a sequence $\{\lambda_k\}_{k\in \mathbb{N}}$ with $\lim_{k\to +\infty}\lambda_k=+\infty$ such that $\lim_{k\to +\infty}v(x_{\lambda_k})< 0$, and, hence, 
\begin{align}\label{eq:minus-infty}
\lim_{k\to +\infty}\dual(\lambda_k)= \lim_{k\to +\infty}
f(x_{\lambda_k})+\lambda_k v(x_{\lambda_k})=-\infty    
\end{align}
since $f(x_{\lambda_k})\le \max_{x\in X}f(x)<\infty$. 

Next, observe that $\dual(0) = \min_{x \in X} f(x)$
is finite because $f(\cdot)$ is continuous and $X$ is compact. Given that the dual function $\dual(\lambda)$, with $\lambda\ge 0$, is continuous because it is concave,  from \eqref{eq:minus-infty} it follows that there exists  $\bar k \in \mathbb{N}$  such that $\dual(\lambda_{\bar k})<\dual(0)$. 
By the Weierstrass' extreme value theorem, the continuous dual function must admit a (global by concavity) maximum over $[0,\lambda_{\bar k}]$.  Moreover, due to the choice of $\lambda_{\bar k}$, this maximum is attained at  some $\ls <\lambda_{\bar k}$, thus implying $\lambda_{\bar k}\notin \Ls$. Since the dual function is concave, then the set $\Ls$  is strictly contained in the interval $[0,\lambda_{\bar k}]$ and hence it is bounded.    This concludes the proof.
\end{pf}
\end{lemma}
\btb

Next, we recall some concepts and results of convex analysis that are useful to motivate the proposed bisection procedure for finding an optimal solution $\ls$ to~\ref{eq:dual_problem} by searching for a zero of the differential mapping of the dual function $\dual(\lambda)$. 

We first need to characterize the differential map of $\dual(\lambda)$. To this end, we recall the following result, which holds since $X$ is compact and $L(x,\lambda)$ is affine in $\lambda$ (cf.~\eqref{eq:lagrangian}), see~\cite[Proposition~B.25]{Bertsekas1997nonlinear}.
\begin{lemma}[Danskin's Theorem] \label{lemma:subdifferential}
The sub-differential set of $\dual(\lambda)$ at a generic $\lambda \ge 0$ is given by
\begin{equation} \label{eq:subg_comp}
\partial \dual(\lambda) = \conv{ \{v(x) :~ x \in X(\lambda)\}},
\end{equation}
where $X(\lambda)$ is defined in~\eqref{eq:primal_recovery}. \hqed
\end{lemma}
As is clear from Lemma~\ref{lemma:subdifferential}, the dual function is not necessarily differentiable at all $\lambda$'s. In particular, whenever the minimization in~\eqref{eq:primal_recovery} has multiple solutions, $\partial \dual(\lambda)$ is not a singleton and $\dual(\lambda)$ is not differentiable at $\lambda \ge 0$.

Since $\dual(\lambda)$ is concave and proper (cf.~\cite[p.~24]{rockafellar1970convex}), we can recall the following result, see~\cite[Theorem~23.5]{rockafellar1970convex}.
\begin{lemma}[Optimality] \label{lemma:optimality}
\rev The dual function \btb $\dual(\lambda)$ achieves its maximum at $\ls \ge 0$ if and only if $\,0 \in \partial \dual(\ls)$. \hqed
\end{lemma}

Finally, leveraging again the fact that the scalar function $\dual(\lambda)$ is concave and proper, by~\cite[p.~238]{rockafellar1970convex}, $\partial \dual(\lambda)$ is monotone and satisfies
\begin{equation} \label{eq:monotonicity}
(v_2 - v_1) (\lambda_2 - \lambda_1) \le 0,
\end{equation}
for all $v_1 \in \partial \dual(\lambda_1)$, $v_2 \in \partial \dual(\lambda_2)$, $\lambda_1 \ge 0$, and $\lambda_2 \ge 0$.

The facts recalled above lead to the following results.

\begin{prop}\label{prop:zero-violation}
If there exist $\tilde\lambda \ge 0$ with $x_{\tilde\lambda}$ satisfying $v(x_{\tilde\lambda}) = 0$, then, $x_{\tilde\lambda}$ is an optimal primal solution, $\tilde\lambda$ is an optimal dual solution, and the duality gap is zero, i.e.,
\begin{equation*}
\max_{\lambda\ge 0} \dual(\lambda)=f^\star.
\end{equation*}
\begin{pf}
By $v(x_{\tilde\lambda}) = 0$, by definition of $x_{\tilde\lambda}$ in~\eqref{eq:primal_recovery} and~\eqref{eq:dual_function}, it follows that 
\begin{equation} \label{eq:phi-lambda}
	f(x_{\tilde\lambda}) = f(x_{\tilde\lambda}) + \tilde \lambda v(x_{\tilde\lambda}) = \phi(\tilde\lambda) \le f(x)+{\tilde\lambda} v(x),
\end{equation}
for all $x \in X$, which implies 
\begin{align*}
	f(x_{\tilde\lambda})\le \min_{\substack{x\in X \\ v(x)\le 0}} \{f(x)+{\tilde\lambda} v(x)\} \le \min_{\substack{x\in X \\ v(x)\le 0}} f(x) = f^\star.
\end{align*}
On the other hand, since $x_{\tilde\lambda} \in \{x\in X: v(x)\le 0\}$, it is feasible for~\ref{eq:problem}, and, thus, 
\begin{align*}
	f(x_{\tilde\lambda}) \ge \min_{\substack{x\in X \\ v(x)\le 0}} f(x) = f^\star. 
\end{align*}
By the last two inequalities and \eqref{eq:phi-lambda}, we have that
\begin{equation*}
	\dual(\tilde\lambda) = f(x_{\tilde\lambda}) = f^\star
\end{equation*}
so that $x_{\tilde\lambda}$ is an optimal primal solution. Moreover, since $\tilde{\lambda} \ge 0$, we have $\max_{\lambda\ge 0} \dual(\lambda) \ge \dual(\tilde\lambda)=f^\star$ and, then, the zero duality gap assertion follows immediately from \eqref{eq:lower}. Also, $\tilde\lambda \in \Ls$ since $\dual(\tilde\lambda)$ achieves the maximum of the dual function. 
\end{pf}
\end{prop}

\begin{prop} \label{prop:violation_sign}
Under Assumption~\ref{ass:strict-feasibility}, we have that 
\begin{subequations} \label{eq:violation_sign}
\begin{align}
	&\hspace{-1mm}\text{if }	\lambda > \ls, \, \forall \ls \in \Ls  \Rightarrow v(x_{\lambda}) < 0, \, \forall x_{\lambda} \in X(\lambda)\label{eq:negative_violation} \\
	&\hspace{-1mm}\text{if }	\lambda < \ls, \, \forall\ls \in \Ls \Rightarrow  v(x_{\lambda}) > 0, \, \forall x_{\lambda} \in X(\lambda)\label{eq:positive_violation}
\end{align}
where $X(\lambda)$ is defined in~\eqref{eq:primal_recovery}. Moreover, if for $\lambda \ge 0$
\begin{align}
	&\exists\, x_{\lambda} \in X(\lambda) :~  v(x_{\lambda}) < 0 \Rightarrow  \lambda \ge \ls, \forall \ls \in \Ls \label{eq:greater_than_ls} \\
	&\exists\, x_{\lambda} \in X(\lambda) :~  v(x_{\lambda}) > 0  \Rightarrow  \lambda \le \ls, \forall \ls \in \Ls. \label{eq:smaller_than_ls}
\end{align}
\end{subequations}
\begin{pf}
Let $\ls$ be an optimal solution to~\ref{eq:dual_problem}, which exists under Assumption~\ref{ass:strict-feasibility} \rev due to Lemma~\ref{lemma:existence} \btb and satisfies $\ls>0$ since we are considering the case when $v(x)\le 0$ is a support constraint. By Lemma~\ref{lemma:optimality} we know that $0 \in \partial\dual (\ls)$. 

Now consider a $\lambda > \ls$ for all $\ls \in \Ls$. Then, setting $\lambda_1 = \ls$ and $\lambda_2 = \lambda$ in~\eqref{eq:monotonicity} and dividing both sides by $\lambda - \ls > 0$ yields $v_2 \le v_1$ for all $v_2 \in \partial\dual(\lambda)$ and $v_1 \in \partial\dual(\ls)$. Setting $v_1 = 0 \in \partial\dual(\ls)$, we then have that any $v_2 \in  \partial\dual(\lambda)$ satisfies $v_2\le 0$ and, owing to Lemma~\ref{lemma:subdifferential},  $v(x_\lambda)\le 0$, for any $x_\lambda \in X(\lambda)$ defined in~\eqref{eq:primal_recovery}. Since $\lambda \notin \Ls$, then, by Lemma \ref{lemma:optimality} $v(x_\lambda)< 0$, for any $x_\lambda \in X(\lambda)$ \eqref{eq:negative_violation}, which concludes the proof of~\eqref{eq:negative_violation}.
Condition~\eqref{eq:positive_violation} can be proven by following similar steps.

Conversely, if $\lambda$ is such that there exists $x_{\lambda} \in X(\lambda)$ with $v(x_{\lambda}) \in \partial\dual(\lambda)$ satisfying $v(x_{\lambda}) < 0$, then setting $\lambda_1 = \ls$, with $\ls \in \Ls$, $\lambda_2 = \lambda$, $v_1 = 0$, and $v_2 = v(x_\lambda)$ in~\eqref{eq:monotonicity} and dividing both sides by $v_2 - v_1 = v(x_\lambda) < 0$ yields $\lambda \ge \ls$, which holds irrespectively of the choice of $\ls \in \Ls$, thus proving~\eqref{eq:greater_than_ls}.
Condition~\eqref{eq:smaller_than_ls} can be proven by following similar steps.
\end{pf}
\end{prop}

We next exploit the monotonicity property \eqref{eq:monotonicity} and Lemmas \ref{lemma:subdifferential} and  \ref{lemma:optimality} to prove Proposition \ref{prop:characterization} characterizing the dual function and the set of its maximizers $\Ls$.  
\begin{prop}\label{prop:characterization}
Under Assumption \ref{ass:strict-feasibility}, the set $\Ls$ is either a singleton or a closed, possibly infinite, interval. Moreover, in the latter case, 
\begin{itemize}
\item[i)] the dual function $\dual: \R \to \R$ 
in \eqref{eq:dual_function} is differentiable in the interior of $\Ls$ and has zero derivative;
\item[ii)] if $\lambda$ is an interior point of $\Ls$, then, any minimizer $x_\lambda \in X(\lambda)$ in \eqref{eq:primal_recovery} is an optimal primal solution 
\item[iii)] if $\lambda$ is a boundary point of $\Ls$, then, there exists $x_\lambda \in X(\lambda)$ which is  an optimal primal solution 
\end{itemize}
\begin{pf}
\rev By Lemma~\ref{lemma:existence}, the set $\Ls$ is non-empty. \btb Suppose that $\Ls$ is not a singleton. Then, it must be a closed interval since the dual function is concave. 
Now let $\lambda>0$ be an interior point of the interval $\Ls$. Then, there exist $\ls_1, \ls_2\in \Ls$ such that $\ls_1<\lambda<\ls_2$ and, by the monotonicity property  \eqref{eq:monotonicity}, any $v \in \partial \dual(\lambda)$ must satisfy 
\begin{align}\label{eq:system}
	\begin{cases}
		v (\lambda - \ls_1) \le 0\\
		v (\lambda - \ls_2) \le 0
	\end{cases}
\end{align}
where we used the fact that $0 \in \partial \dual(\ls_1) \cap \partial \dual(\ls_2)$ (cf. Lemma \ref{lemma:optimality}). 
Since $\ls_1 < \lambda < \ls_2$, then, $0$ is the only possible value of $v$ satisfying~\eqref{eq:system} and, hence, $\partial \dual(\lambda)=\{0\}$. This concludes the proof of the differentiability result i). \\
If we recall the characterization of the sub-differential set in \eqref{eq:subg_comp}, we then have that for any interior point $\lambda$ of $\Ls$, $v(x_\lambda)=0$ for all $x_\lambda \in X(\lambda)$,  and point ii) immediately follows from Proposition \ref{prop:zero-violation}. If $\lambda$ is the left extreme of $\Ls$ and $\ls$ is an interior point, by the monotonicity property any $v \in \partial \dual(\lambda)$ must satisfy
\begin{equation*}
	v (\lambda-\ls)\le 0,
\end{equation*}
which implies that $ v \ge 0$. A similar reasoning applies to the right extreme of $\Ls$ (if $\Ls$ is bounded) for which $v\le 0$. Since in both cases $0 \in \partial \dual(\lambda)$ (cf. Lemma \ref{lemma:optimality}), then, by Lemma \ref{lemma:subdifferential} we have that there exists $x_\lambda$ with $v(x_\lambda) = 0$, which, by Proposition \ref{prop:zero-violation}, is an optimal primal solution. This concludes the discussion of point iii) and the proof.
\end{pf}
\end{prop}

Thanks to the monotonicity property in~\eqref{eq:monotonicity}, we can also state the following result about the cost achieved by any primal solution computed via~\eqref{eq:primal_recovery}.
\begin{prop} \label{prop:cost_monotonicity}
For any $\lambda_1 \ge 0$ and $\lambda_2 \ge 0$ we have
\begin{equation} \label{eq:cost_monotonicity}
\lambda_1 < \lambda_2 \Rightarrow f(x_{\lambda_1}) \le f(x_{\lambda_2}) ~\land~ v(x_{\lambda_1}) \ge v(x_{\lambda_2}),
\end{equation}
for all $x_{\lambda} \in X(\lambda)$ defined in~\eqref{eq:primal_recovery}.
\begin{pf}
By~\eqref{eq:primal_recovery} and the definition of minimizer, any $x_{\lambda_1}\in X(\lambda_1)$ satisfies
\begin{align}
	f(x_{\lambda_1}) + \lambda_1 v(x_{\lambda_1})
	&\le 	f(x_{\lambda_2}) + \lambda_1 v(x_{\lambda_2}), \label{eq:min_cost_inequality}
\end{align}
for any $x_{\lambda_2} \in X(\lambda_2) \subseteq X$.

If $\lambda_1 < \lambda_2$, then, dividing both sides of~\eqref{eq:monotonicity} by $\lambda_2 - \lambda_1 > 0$ yields $v_2 - v_1 \le 0$ for all $v_2 \in \partial\dual(\lambda_2)$ and $v_1 \in \partial\dual(\lambda_1)$. Owing to Lemma~\ref{lemma:subdifferential} and~\eqref{eq:primal_recovery}, we can select $v_1 = v(x_{\lambda_1})$ and $v_2 = v(x_{\lambda_2})$, thus obtaining
\begin{equation} \label{eq:violation_inequality}
	v(x_{\lambda_2}) - v(x_{\lambda_1}) \le 0.
\end{equation}
Rearranging the terms in~\eqref{eq:min_cost_inequality} we have
\begin{equation*}
	f(x_{\lambda_1}) - f(x_{\lambda_2}) \le \lambda_1 (v(x_{\lambda_2}) - v(x_{\lambda_1})) \le 0,
\end{equation*}
the second inequality being due to~\eqref{eq:violation_inequality} and $\lambda_1 \ge 0$, thus concluding the proof.
\end{pf} 
\end{prop}

\subsection{Proof of Theorem~\ref{thm:main}} \label{sec:proof}
We start by showing that the first cycle in Algorithm~\ref{alg:dual-bisection} returns a feasible solution after a finite number of iterations $K \ge 0$.

If $\lambdainit$ is such that $v(x_{\lup(0)}) = v(x_{\lambdainit}) \le 0$, then the first cycle (Steps~\ref{step:cycle_1_start}-\ref{step:cycle_1_end}) is skipped and the \rev first part of the \btb statement is trivially satisfied with $K = 0$. Moreover, under the standing assumption that $v(x_0) > 0$ we have that $v(x_{\ldn(0)}) = v(x_0) > 0$. Otherwise, if $v(x_{\lambdainit}) > 0$, we enter the first cycle at least once ($K > 0$) and keep updating $\lup$ and $\ldn$ according to the following dynamical system
\begin{subequations} \label{eq:lupk_ldnk_explicit}
\begin{align}
&\lup(\iter+1) = 2 \, \lup(\iter) = \cdots = 2^{\iter+1} \lup(0), \label{eq:lupk_explicit} \\
&\ldn(\iter+1) = \lup(\iter) = 2^\iter \lup(0),  \label{eq:ldnk_explicit}
\end{align}
\end{subequations}
with $\lup(0) = \lambdainit$ and $\ldn(0) = 0$, as long as $v(x_{\lup(\iter)}) > 0$. Since $\ldn(\iter+1) = \lup(\iter)$, then $v(x_{\ldn(\iter+1)}) = v(x_{\lup(\iter)}) > 0$ and, according to~\eqref{eq:smaller_than_ls} of Proposition~\ref{prop:violation_sign},
\begin{equation} \label{eq:bound_on_ldnkp}
\hspace{-2.5mm} v(x_{\lup(\iter)}) > 0 \Rightarrow \lup(\iter) = \ldn(\iter+1) \le \ls \,,  \forall \ls \in \Ls.
\end{equation}
Since, by Assumption~\ref{ass:strict-feasibility} \rev and Lemma~\ref{lemma:existence}, \btb $\Ls \subset [0,+\infty)$ is non-empty, and is closed as an effect of $\dual(\lambda)$ being a concave function, it admits a minimum $\lsdn = \min_{\ls \in \Ls} \ls$. Using~\eqref{eq:ldnk_explicit} in~\eqref{eq:bound_on_ldnkp} and recalling $\lup(0) = \lambdainit$ we obtain
\begin{equation} \label{eq:cycling_condition}
v(x_{\lup(\iter)}) > 0 \Rightarrow 2^\iter \lambdainit \le \lsdn \Leftrightarrow \iter \le \log_2 \frac{\lsdn}{\lambdainit},
\end{equation}
which entails that the first cycle terminates after $K = \lfloor \log_2 \frac{\lsdn}{\lambdainit} \rfloor + 1$ iterations \rev with 
$v(x_{\ldn(K)}) > 0$. \btb Moreover, since $K > \log_2 \frac{\lsdn}{\lambdainit}$, then $\lup(K) = 2^{K} \lambdainit > \lsdn$ and, by reversing~\eqref{eq:cycling_condition}, we have that $\lup(K)$ must satisfy $v(x_{\lup(K)}) \le 0$, meaning that $x_{\lup(K)}$ is a feasible solution for~\ref{eq:problem}.

The previous discussion shows also that the first element $\hat{x}(K) = x_{\lup(K)}$ of the sequence $\{\hat{x}(\iter)\}_{\iter\ge K}$ is feasible for~\ref{eq:problem}. Moreover, if $\lup(K)$ happens to \rev belong to \btb the interior of $\Ls$, then by Proposition~\ref{prop:characterization}.ii $x_{\lup(K)}$ \rev is \btb optimal for~\ref{eq:problem} and \rev is \btb returned by Algorithm~\ref{alg:dual-bisection} since $v(x_{\lup(K)}) = 0$. If $\lup(K)$\rev, instead, \btb happens to land on the right extreme of $\Ls$, then by Proposition~\ref{prop:characterization}.iii $x_{\lup(K)}$ may be optimal for~\ref{eq:problem} and is returned by Algorithm~\ref{alg:dual-bisection} if $v(x_{\lup(K)}) = 0$. In both cases, we can easily see that the second part of the statement of Theorem~\ref{thm:main} holds with $\bar{k} = K$.

In all other cases we have $v(x_{\lup(K)}) < 0$ and the execution proceeds with the bisection cycle (Steps~\ref{step:cycle_2_start}-\ref{step:cycle_2_end}). Since $\hat{x}(K)$ is feasible and $\hat{x}(\iter+1)$ is set equal to $\hat{x}(\iter)$ or to a $x_{\lambda(\iter)} = x_{\lup(\iter)}$ if $v(x_{\lambda(\iter)}) < 0$, then the sequence $\{\hat{x}(\iter)\}_{\iter \ge K}$ satisfies~\eqref{eq:feasible} by construction. Again by construction, for all $\iter \ge K$, $\lambda(\iter) < \lup(\iter)$ and, by Steps~\ref{step:lup_update} and~\ref{step:lup_keep}, it is easy to see that either $\lup(\iter+1) < \lup(\iter)$ and, by Proposition~\ref{prop:cost_monotonicity}, $f(\hat{x}(\iter+1)) \le f(\hat{x}(\iter))$, or $\lup(\iter+1) = \lup(\iter)$ and $\hat{x}(\iter+1) = \hat{x}(\iter)$, and thus $f(\hat{x}(\iter+1)) = f(\hat{x}(\iter))$, which shows~\eqref{eq:decreasing_cost}.

To conclude the proof we shall note that $\lambda(\iter)$ may land inside $\Ls$, thus yielding an optimal solution $x_{\lambda(\iter)}$ with \rev $v(x_{\lambda(\iter)})= 0$, \btb or may land on the boundary of $\Ls$, which \emph{may} yield an optimal solution if $v(x_{\lambda(\iter)}) = 0$ (cf. Proposition~\ref{prop:characterization}). In these cases, the algorithm is halted, and $\bar{k}$ in the statement of Theorem~\ref{thm:main} is equal to the current iteration index. This last observation concludes the proof. \hqed

\subsection{Proof of Theorem~\ref{thm:skip_first_cycle}} \label{sec:proof_skip_first_cycle}

If a solution $\slater \in X$ such that $v(\slater) < 0$ is known, then, for all $\ls \in \Ls$, we have
\begin{align} 
    \dual(0) < \dual(\ls) &= \min_{x \in X} f(x) + \ls v(x) \notag \\
    &\le f(\slater) + \ls v(\slater), \label{eq:slater_ineq}
\end{align}
where the first inequality is due to $0 \not\in \Ls$, the equality is due to~\eqref{eq:lagrangian} and~\eqref{eq:dual_function}, and the last inequality holds by definition of minimum.
Relation~\eqref{eq:slater_ineq} implies 
\begin{equation*}
\ls < \frac{\dual(0)-f(\slater)}{v(\slater)} = \lambdainit,
\end{equation*}
for all $\ls \in \Ls$, the equality being due to~\eqref{eq:slater_init}.

By~\eqref{eq:negative_violation} in Proposition~\ref{prop:violation_sign}, $v(x_{\lup(0)}) = v(x_{\lambdainit}) < 0$ and the first cycle is skipped, thus yielding $K = 0$ and concluding the proof. \hqed

\section{Theoretical analysis of the algorithm in \cite{Manieri-IFAC-2023}}
\label{app:multi-agent}
\subsection{Proof of Proposition~\ref{prop:ifac-inner-loop-scalar}}
 For the sake of contradiction, suppose that $\exists \breve\kappa\ge 0$ such that for all $\kappa\ge\breve \kappa$ the solutions $\xi_1(\kappa),\dots,\xi_m(\kappa)$ computed in Step~\ref{step:XLP_tentativePrimal} of Algorithm~\ref{algo:ifac-inner} jointly violate the coupling constraint~\eqref{eq:multi-agent MILP coupling constraints}, i.e. $\sumim A_i \xi_i(\kappa)  -b > 0$.

Denote with $\mathcal{V}$ the set of all positive values that the scalar violation $\sumim A_i \xi_i  -b $ can assume as $\xi_i \in\vertt{\Xi_i}$ for $\im$, i.e.
    \begin{equation}
        \label{eq:ifac-violation-set}
        \begin{aligned}
             \cV = \Big\{ v\in \R :~ &v = \sumim A_i \xi_i  -b \wedge v>0 , \\ &\xi_i \in \vertt{\Xi_i} \,\forall \, \im\Big\} ,
        \end{aligned}
    \end{equation}
    which is finite due to $\Xi_i, \im$ being bounded. Let $\tunderbar{v}=\min_{v\in\cV} v > 0$ be the minimum (scalar) value in $\cV$, then for $\kappa \ge \breve \kappa$ we can lower-bound the expression of $\lambda(\kappa)$ in Step~\ref{step:ifac-dual-update} of Algorithm~\ref{algo:ifac-inner} as
    \begin{align}
        \label{eq:proof-lambda-LB}
        \lambda(\kappa) &\ge \lambda(\breve \kappa)+\sum_{\tau=\breve \kappa}^{\kappa-1}\alpha(\tau) \Big ( \sumim A_i \xi_i(\tau) -b\Big) \notag \\
        &\ge \sum_{\tau=\breve \kappa}^{\kappa-1}\alpha(\tau) \tunderbar{v} =  \tunderbar{v} \sum_{t=\breve \kappa}^{\kappa-1}\alpha(\tau)\, ,
    \end{align} 
    where the first inequality is due to the definition of the projection operator as $[\cdot]_+=\max\{0,\cdot\}$, whereas the second inequality follows from $\lambda(\breve \kappa)\ge0$ and the assumption by contradiction that $\sumim A_i \xi_i
    (\kappa) -b>0$ for $\kappa>\breve \kappa$ jointly with the definition of $\tunderbar{v}$. 

    Since the step-size $\alpha(\kappa)$ satisfies~\eqref{eq:alphak}, then \eqref{eq:proof-lambda-LB} implies that $\lim_{\kappa\to\infty} \lambda(\kappa) =\infty$. This, however, contradicts Theorem 2.2~in~\cite{shor2012minimization}, which guarantees convergence of the sequence $\{\lambda(\kappa)\}_{\kappa\ge0}$ when~\eqref{eq:alphak} holds. Thus, we have that
    \begin{equation*}
        \forall \breve{\kappa}\ge 0 \,\,\, \exists \bar{\kappa}\ge\breve{\kappa} : \sumim A_i \xi_i(\bar \kappa) - b \le 0 \, ,
    \end{equation*}
    which concludes the proof. \hqed

\end{document}